\documentclass[a4paper,12pt]{article}

\usepackage{amssymb,amsmath,amsfonts, amsthm}
\usepackage{graphicx}

\newcommand{\R}{\mathbb{R}}
\newcommand{\N}{\mathbb{N}}
\newcommand{\Z}{\mathbb{Z}}
\newcommand{\Q}{\mathbb{Q}}
\newcommand{\C}{\mathbb{C}}
\newcommand{\bbH}{\mathbb{H}}

\newcommand{\tr}{\mathrm{tr}}
\newcommand{\Tr}{\mathrm{Tr}}

\newtheorem{theorem}{Theorem}
\newtheorem{proposition}[theorem]{Proposition}
\newtheorem{corollary}[theorem]{Corollary}
\newtheorem{lemma}[theorem]{Lemma}
\newtheorem{conjecture}[theorem]{Conjecture}
\newtheorem{claim}[theorem]{Claim}

\theoremstyle{definition}

\theoremstyle{remark}
\newtheorem*{rem}{Remark}

\title{A GEOMETRIC CHARACTERIZATION OF ARITHMETIC FUCHSIAN GROUPS}
\author{Slavyana Geninska, Enrico Leuzinger}

\begin{document}

\maketitle
\begin{abstract}
The trace set of a Fuchsian group $\Gamma$  ist the set of length of closed geodesics in the surface $\Gamma\backslash \bbH$.  Luo and Sarnak showed that the trace set of a cofinite arithmetic Fuchsian group satisfies the bounded clustering property. Sarnak then conjectured that the B-C property actually characterizes  arithmetic Fuchsian groups. 
Schmutz stated the even stronger conjecture that a cofinite Fuchsian group  is arithmetic if its trace set  has linear growth. He proposed a proof of this conjecture in the case when the group $\Gamma$ contains at least one parabolic element, 
but unfortunately this proof contains a gap.
In the present paper we point out this  gap and  we prove Sarnak's conjecture under 
the assumption that the Fuchsian group $\Gamma$ contains  parabolic elements.  
\end{abstract}

\section{Introduction}
\label{intro}

Let $\Gamma$ be a  Fuchsian group, i.e. a discrete subgroup of $PSL(2, \mathbb R)$.  Such a $\Gamma$ acts properly discontinuously and  isometrically on the hyperbolic plane $\bbH$ and $M = \Gamma \backslash \bbH$ is a Riemann surface. The \emph{trace set of $\Gamma$} and the \emph{trace set of $M$} are defined as follows:
\begin{eqnarray*}\Tr(\Gamma) & := & \{ \tr(T) \mid T\in \Gamma \},\\
\Tr(M) & = & \{ 2\cosh \frac{L(a)}{2} \mid a~ \mathrm{is~ a~ closed~ geodesic~in} ~M~ \mathrm{of~length} ~ L(a)\}. 
\end{eqnarray*}
These two subsets of $\R$ in fact coincide for torsionfree $\Gamma$ (see Section 2.1 below).

It is a general question  if certain classes of Fuchsian groups can be characterized by means  of their trace set or, equivalently, by the trace set of the surfaces that they define.
In this paper we are interested in  characterizations of arithmetic Fuchsian groups. 
There is a classical  characterization of arithmetic Fuchsian groups
due to Takeuchi which is based on number theoretical properties of  their trace sets  \cite{kT75}. 

 Luo and Sarnak pointed out large scale properties of the behaviour of the trace set of arithmetic Fuchsian groups.
We say that the trace set of a Fuchsian group $\Gamma$ satisfies the \textit{bounded clustering} or \textit{B-C property} iff there exists a constant $B(\Gamma)$ such that for all integers $n$ the set $\Tr(\Gamma) \cap [n, n+1]$ has less than $B(\Gamma)$ elements. Further set
$$ Gap(\Gamma) := inf \{ |a-b| \mid a,b \in \Tr(\Gamma), a\neq b \}. $$

In \cite{wL94} Luo and Sarnak made a first step towards a new geometric characterization of arithmetic Fuchsian groups
 by proving the following result:

\begin{theorem}[\cite{wL94}]
\label{T:LuoSarnak}
Let $\Gamma$ be a cofinite Fuchsian group.
\begin{itemize}
\item[(i)] If $\Gamma$ is  arithmetic then $\Tr(\Gamma)$ satisfies the B-C property.
\item[(ii)] If $\Gamma$ is  derived from a quaternion algebra then $Gap(\Gamma)>0$.
\end{itemize}
\end{theorem}
 
 Sarnak  conjectured  that the converse of Theorem \ref{T:LuoSarnak} also holds. 

\begin{conjecture}[Sarnak \cite{pS95}] Let $\Gamma$ be a cofinite Fuchsian group.
\label{Con:Sarnak}
\begin{itemize}
\item[(i)] If $\Tr(\Gamma)$ satisfies the B-C property then $\Gamma$ is  arithmetic. 
\item[(ii)] If $Gap(\Gamma) > 0$ then $\Gamma$ is  derived from a quaternion algebra.
\end{itemize}
\end{conjecture}

In \cite{pS96}  Schmutz makes an even stronger conjecture using  the linear growth of a trace set instead of the B-C property.
The trace set of a Fuchsian group $\Gamma$ is said to have \textit{linear growth} iff there exist positive real constants $C$ and $ D$
 such that for every $n \in \N$
$$ \# \{ a \in \Tr(\Gamma) \mid a \leq n \} \leq D + nC.$$ 

\begin{rem} If a Fuchsian group $\Tr(\Gamma)$ satisfies the B-C property, then $\Tr(\Gamma)$ has linear growth with $D = 0$ and $ C = B(\Gamma)$.
But the opposite is not true in general: B-C $\nRightarrow$ linear.
\end{rem}

\begin{conjecture}[Schmutz \cite{pS96}] 
\label{Con:Schmutz}
Let $\Gamma$ be a cofinite Fuchsian group. If $\Tr(\Gamma)$ has linear growth then $\Gamma$ is  arithmetic.
\end{conjecture}

In \cite{pS96} Schmutz proposed a proof of Conjecture \ref{Con:Schmutz} in the case when $\Gamma$ contains at least one parabolic element.
But unfortunately the proof contains a gap as we will point out  in  Section \ref{growth}.

It remains an open question wether the gap in \cite{pS96} can be closed. Observe that a positive answer would imply that there
are cofinite Fuchsian groups (with parabolic elements) whose trace set grows linearly but without satisfying the stronger B-C property.
Furthermore we remark the conjectures of Sarnak and Schmutz remain completely open for cocompact Fuchsian groups.

\paragraph{ }
The plan of the paper is as follows. 
In Section \ref{defs} we fix the notation and give some basic definitions and results.
In Section \ref{ypieces} we prove (or list) some auxilary results that are used later.
In the last Section \ref{growth} we use techniques similar to those developed by Schmutz to prove part (a) of Sarnak's conjecture
under the assumption  that the Fuchsian group $\Gamma$ contains at least one parabolic element.

\section{Some basic definitions and facts}
\label{defs}
\setcounter{theorem}{0}

\subsection{Trace sets}
A general reference for this section is  the book  \cite{sK92}.
We denote by $SL(2,\R)$ the group of real $2 \times 2$ matrices with determinant $1$ and by $PSL(2, \R)$ the
 quotient group $SL(2, \R))/\{ \pm 1_2\}$ where $1_2$ is the $2 \times 2$ identity matrix.

A \textit{Fuchsian group} is a discrete subgroup of $PSL(2,\R)$. On the 
 \textit{hyperbolic plane} $\bbH = \{ z = x + iy \in \C \mid y > 0 \}$ endowded with the metric $ds² = y^{-2}(dx^2 +dy^2)$ 
a Fuchsian group  acts isometrically and properly discontinuously by fractional linear transformations
\[ \{z \mapsto \frac{az+b}{cz+d} \mid a,b,c,d \in \R, ~ad-bc=1\} \subseteq Isom(\bbH).\]
\paragraph{ }
For $T = \begin{bmatrix} a&b\\c&d \end{bmatrix} \in PSL(2,\R)$ we set $\tr(T):=|a+d|$.  
For a Fuchsian group $\Gamma$ we then call 
\[ \Tr(\Gamma) = \{\tr(T) \mid T \in \Gamma \} \]
the \textit{trace set of} $\Gamma$.
\paragraph{ }
Let $T = \begin{bmatrix} a&b\\c&d \end{bmatrix} \in PSL(2,\R)$ with $c \neq 0$. The circle
\[ I(T) = \{z \in \C \mid |cz+d| = 1 \}, \]
which is the subset of $\C$, where $T$ acts as an Euclidean isometry, is called the \textit{isometric circle} of $T$. It is uniquely determined by $c$ and $d$, because its center is $(-\frac{d}{c},0)$ and its radius is equal to $\frac{1}{|c|}$. 

\begin{theorem}
Let $T = \begin{bmatrix} a&b\\c&d \end{bmatrix} \in PSL(2,\R)$ with $c \neq 0$. The isometric circles $I(T)$ and $I(T^{-1})$ have the same radius; and $I(T^{-1}) = T(I(T))$.It can also be used as a definition for cofinite arithmetic Fuchsian groups.
\end{theorem}

\paragraph{ }
A Fuchsian group $\Gamma$  is called \textit{cofinite} or \textit{of the first kind} if the associated quotient surface has finite area,  $\mu(\Gamma \backslash \bbH) < \infty$.

\paragraph{ }
For a Fuchsian group $\Gamma$ let $M$ be the quotient $\Gamma \backslash \bbH$ with the points corresponding to fixed points of elliptic elements in $\Gamma$ removed. Endowed  with the metric induced by the hyperbolic metric on $\bbH$,  $M$ is a Riemann surface.

\paragraph{ }
Let $a$ be a closed geodesic on $M$. Then, by abuse of notation, the length of $a$ is also denoted by $a$.
We define the \textit{trace of} $a$ to be $\tr(a):=2 \cosh \frac{a}{2}$ and we set  $\Tr(M)=\{ \tr(a) \mid a~ \mathrm{is~ a~ closed~ geodesic~in} ~ M\}$.

\begin{proposition} 
For a torsionfree Fuchsian group $\Gamma$ holds $\Tr(\Gamma)=\Tr(\Gamma \backslash\bbH)$.
\end{proposition}
\begin{proof} 
Let  $C(T)$ be the axis  of a hyperbolic element $T$ in $\Gamma$.
The image of $C(T)$ in $M=\Gamma \backslash \bbH$   is a closed geodesic $a$  with length equal to the distance between $x$ and $T(x)$ for any $x \in C(T)$. Vice versa, for every closed geodesic $a$ on $M$ there exists a hyperbolic element $T$ in $\Gamma$, such that the image in $M$ of its axis is $a$, and for every $x \in C(T)$ the distance between $x$ and $T(x)$ is $a$.

Every hyperbolic element $T$ in $\Gamma$ can be conjugated by a hyperbolic isometry $R$ to $T' = RTR^{-1} = \begin{bmatrix} e ^{\frac{\tau}{2}}& 0 \\ 0 & e ^{-\frac{\tau}{2}}\end{bmatrix}$ for some $\tau > 0$. Then we have on the one hand 
\[ \tr(T)=\tr(T') =  e ^{\frac{\tau}{2}} +e ^{-\frac{\tau}{2}} = 2\cosh\frac{\tau}{2}.\] 
On the other hand $T(R^{-1}(i))= R^{-1}(T'(i))=R^{-1}(e^\tau i)$. Hence, as $R^{-1}(i) \in C(T)$, the length of the closed geodesic on $M$ defined by $T$ is equal to 
\[d_{h}(R^{-1}(i), T(R^{-1}(i))) = d_h(i, e^{\tau} i) = \ln \frac{|e^\tau i+i| + |e^\tau i-i|}{|e^\tau i+i| - |e^\tau i-i|} = \ln\frac{2e^\tau}{2}= \tau.\]
\end{proof}

\subsection{Takeuchi's characterization of arithmetic Fuchsian groups}
In order to state Takeuchi's results we recall some definitons and facts concerning quaternion algebras.
For more details  we refer to \cite{sK92}, Chapter 5, and to \cite{cM03}, Chapter 0. In this section $F$ will always denote a general field.
\paragraph{ }
Recall that a \textit{quaternion algebra over $F$} is a central simple algebra over $F$ which is four dimensional $F$-vector space.
Each quaternion algebra is isomorphic to an algebra $A=\left( \frac{a,b}{F} \right)$ with $a,b \in F^*=F-\{0\}$ and a basis $\{1,i,j,k\}$, where $i^2 = a$, $j^2=b$, $k = ij = -ji$. 
If each element of a quaternion algebra $A$ has an inverse, then A is called a \textit{division quaternion algebra}.

\paragraph{ }
If $F$ is an algebraic number field it can be written as $\Q(t)$, where $t$ satisfies a polynomial with rational coefficients and $\Q(t)$ is the smallest field containing $\Q$ and $t$. Let $f \in \Q[x]$ be the minimal polynomial of $t$. If $n$ is the dimension of $F$ considered as a vector space over $\Q$, then $f$ has degree $n$.
Let $t_1 = t$, $t_2$, ... , $t_n$ denote the roots of $f$, then the substitution $t \to t_i$ induces a field isomorphism $\Q(t) \to \Q(t_i)$. Conversely, if $\sigma: F = \Q(t) \to \C$ is a field \textit{embedding}, i.e. $\sigma: F \to \sigma(F)$ is a field isomorphism, then $\sigma(t)$ is a root of the minimal polynomial of $t$. Therefore, there are exactly $n$ field embeddings $\sigma: F \to \C$.
$F$ is a \textit{totally real algebraic number field} iff for each embedding of $F$ into $\C$ the image lies inside $\R$.
An element of $F$ is an \textsl{algebraic integer} iff it satisfies a polynomial with coefficients in $F$ and leading coefficient $1$. The algebraic integers of $F$ form a ring and we denote it by $\mathcal{O}_F$.

\paragraph{ }
Let $A=\left( \frac{a,b}{F} \right)$ be a quaternion algebra. For every $x \in A$, $x = x_0+x_1i+x_2j+x_3k$, we define the \textit{reduced norm} of $x$ to be $Nrd(x) = x \bar{x} = x_0^2 - x_1^2 a - x_2^2 b + x_3^2 ab$, where $\bar{x}= x_0-x_1i-x_2j-x_3k$.
An \textit{order} $\mathcal{O}$ in a quaternion algebra $A$ over $F$ is a subring of $A$ containing $1$, which is a finitely generated $\mathcal{O}_F$-module and generates the algebra $A$ over $F$.
The \textit{group of units} in $\mathcal{O}$ of reduced norm $1$ is $\mathcal{O}^1 = \{ x \in \mathcal{O} \mid Nrd(x) = 1\}$.

\paragraph{ }
\label{D:GammaAO}
Let $F$ be a totally real algebraic number field of degree $n$ and let $\varphi_i$, $i \in \{1 \ldots n \}$, be the $n$ distinct embeddings of $F$ into $\C$, where $\varphi_1 = id$. Let $A=\left( \frac{a,b}{F} \right)$ be a quaternion algebra over $F$ such that there exist $n$ $\R$-isomorphisms: for $2 \le i \le n$
$$\rho_i: \left( \frac{\varphi_i(a),\varphi_i(b)}{\R} \right) \rightarrow \left( \frac{-1,-1}{\R} \right) ~ \mathrm{and} ~ \rho_1: \left( \frac{\varphi_1(a),\varphi_1(b)}{\R} \right) \rightarrow M(2, \R).$$

\begin{theorem}
\label{T:G(A,O)FuchsianGroup}
$\Gamma(A,\mathcal{O}):=\rho_1(\mathcal{O}^1)/\{+1_2, -1_2\}$ is a Fuchsian group.
\end{theorem}

 A Fuchsian group $\Gamma$ is \textit{derived from a quaternion algebra} iff $\Gamma$ is a subgroup of finite index of some $\Gamma(A,\mathcal{O})$.
Two Fuchsian groups are \textit{commensurable} iff their intersection has finite index in each of them.
A Fuchsian group 
$\Gamma$ is  \textit{arithmetic} iff $\Gamma$ is commensurable with some $\Gamma(A,\mathcal{O})$.

\paragraph{ }
The following two theorems due to Takeuchi provide an algebraic characterization of (cofinite) arithmetic Fuchsian groups. 

\begin{theorem}[\cite{kT75}]
\label{T:TakeuchiArithmQuat}
Let $\Gamma$ be a cofinite Fuchsian group. Then $\Gamma$ is derived from a quaternion algebra over a totally real algebraic number field if and only if $\Gamma$ satisfies the following two conditions:
\begin{itemize}
\item[(i)] $K:= \Q(\Tr(\Gamma))$ is an algebraic number field of finite degree and $\Tr(\Gamma)$ is contained in the ring of integers $\mathcal{O}_K$ of $K$.
\item[(ii)] For any embedding $\varphi$ of $K$ into $\C$, which is not the identity, $\varphi(\Tr(\Gamma))$ is bounded in $\C$.
\end{itemize}
\end{theorem}

\begin{theorem}[\cite{kT75}]
\label{T:TakeuchiArithm}
Let $\Gamma$ be a cofinite Fuchsian group and $\Gamma^{(2)}$ be the subgroup of $\Gamma$ generated by the set $\{T^2 \mid T \in \Gamma\}$. Then $\Gamma$ is an arithmetic Fuchsian group if and only if $\Gamma^{(2)}$ is derived from a quaternion algebra.
\end{theorem}


\section{Y-pieces and lengths of geodesics on them}
\label{ypieces}
\setcounter{theorem}{0}

\paragraph{ }
	An \textit{Y-piece} is a surface of signature (0, 3), i.e. homeomorphic to a topological sphere with three points removed. For non-negative real numbers 
	$a$, $b$, $c$ we denote with $Y(a, b, c)$ an Y-piece with boundary geodesics of lengths $a$, $b$, $c$.
It is well known that for given boundary geodesics the Y-piece is uniquely determined up to isometry, see \cite{pB92}, Theorem 3.1.7.
We will say that an Y-piece $Y(u,v,w)$ is contained in another Y-piece $Y(x,y,z)$ iff $Tr(Y(u,v,w)) \subseteq Tr(Y(x,y,z))$.

\subsection{Generation of Y-pieces}
\label{SS:Genypieces}
	In this section we are going to show that for every Y-piece one
	 can find a Fuchsian group $\Gamma$ generated by only two elements such that $\Gamma \backslash \bbH$ contains the Y-piece.
	
	The next two lemmas are needed as a preparation for the proof of Proposition \ref{P:generationParab}.
	\begin{lemma} \label{L:CIorthogonal}
	For every hyperbolic transformation $T = 
	    \begin{bmatrix}
			 a & b \\
			 c & d
			\end{bmatrix}$, $c\neq0$, the axis $C(T)$ and the isometric circle $I(T)$ intersect orthogonally. 
	\end{lemma}
	\begin{proof}
	The isometric circle $I(T)$ is the circle with center $M:=(-\frac{d}{c}, 0)$ and radius $r_1 := \frac{1}{|c|}$. 
	The axis $C(T)$ considered as an Euclidean circle has center $N := (\frac{a-d}{2c}, 0)$ and radius $r_2 := \sqrt{\frac{(a+d)^2 - 4}{4c^2}}$. 
	The distance between $M$ and $N$ is $\left| \frac{a+d}{2c} \right|$. 
	 
	 Since $|MN|>r_1$, $|MN|>r_2$ and $|MN|< r_1 + r_2$,  there exists an Euclidean triangle with sides equal to $r_1$, $r_2$ and $|MN|$. Hence $I(T)$ and $C(T)$ intersect. Let $P$ be their intersection point. Then the Euclidean triangle $MNP$ has a right angle at $P$. Indeed:
			$|PN|^2 = \frac{(a+d)^2 - 4}{4c^2}, \;
			  |PM|^2 = \frac{1}{c^2}\  \textup{and}\  \;
			  |MN|^2 = \left( \frac{a-d}{2c} + \frac{2d}{2c} \right)^2 =
			           \frac{(a+d)^2}{4c^2} $.
	\end{proof}
  \begin{lemma} \label{L:Notintersect}
     Let $T_1 = \begin{bmatrix} a_1 & b_1 \\ c & d \end{bmatrix}$ and 
     $T_2 = \begin{bmatrix} a_2 & b_2 \\ \lambda c & \lambda d \end{bmatrix}$ be hyperbolic 
      isometries in $PSL(2, \R)$, where $\lambda$ is a positive real number and $c \neq 0$. If $(a_1 +d)(a_2 + \lambda d) < 0$ then the axes $C(T_1)$ and $C(T_2)$ do not intersect and  on the real axis the repulsive fixed points of $T_1$ and $T_2$ are between their attracting fixed points.
  \end{lemma}
  \begin{proof}
  
  \begin{figure}
  \centering
  \includegraphics[width=1.0\textwidth]{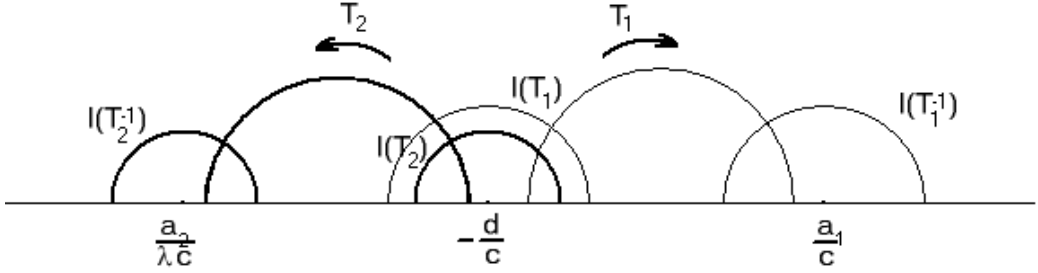}
   \caption{The case $\frac{a_1+d}{c}>0$.}
   \label{F:Notintersect}
  \end{figure}
  
  The isometric circles $I(T_1)$ and $I(T_2)$ have the same center $(-\frac{d}{c}, 0)$. The center of $I(T_1^{-1})$ is $(\frac{a_1}{c}, 0)$ and the center of $I(T_2^{-1})$ is $(\frac{a_2}{\lambda c}, 0)$.
   If $\frac{a_1+d}{c}>0$ then $\frac{a_2+\lambda d}{c}<0$ and hence $\frac{a_2}{\lambda c}<-\frac{d}{c}<\frac{a_1}{c}$. Analogously if $\frac{a_1+d}{c}<0$ then $\frac{a_1}{c}<-\frac{d}{c}<\frac{a_2}{\lambda c}$ (Fig \ref{F:Notintersect}).
  
By Lemma \ref{L:CIorthogonal}  $C(T_1)$ is orthogonal to $I(T_1)$. Hence the radius of $C(T_1)$ considered as an Euclidean circle is shorter than the distance between the centers of $C(T_1)$ and $I(T_1)$. Therefore the repulsive fixed point of $T_1$ lies between the centers of $I(T_1)$ and $C(T_1)$ (Fig \ref{F:Notintersect}). Similarly the attracting fixed point of $T_1$ (which is also the repulsive fixed point of $T_1^{-1}$) lies between the centers of $I(T_1^{-1})$ and $C(T_1)$. Analogous considerations for $T_2$ prove the lemma.
  \end{proof}
 In the next Proposition \ref{P:generationParab} we give sufficient conditions for two hyperbolic isometries to generate a group $\Gamma$ such that $\Gamma \backslash \bbH$ contains an Y-piece $Y(u,v,0)$.

	\begin{proposition} \label{P:generationParab}
	Let $u$ and $v$ be non-negative real numbers. Further let $T_u$ and $T_v$ be elements of $PSL(2, \R)$ such that
	\(
		T_{u} = \begin{bmatrix}
			 a_1 & b_1 \\
			 c & d
			\end{bmatrix}
	\) 
	and
	\(
		T_{v} = 
			\begin{bmatrix}
			 a_2 & b_{2} \\
			 c & d
			\end{bmatrix}
	\),
	with $c \neq 0$ and such that  for $\varepsilon = \pm 1$, 
	$a_1 + d = \varepsilon \tr(u)$ and $a_2 + d = - \varepsilon \tr(v)$.
	Then $\Gamma = \langle T_{u},T_{v} \rangle$ is a Fuchsian group and the surface $\Gamma \backslash \bbH$ contains an Y-piece $Y(u,v,0)$. 
	\end{proposition}
	
	\begin{proof}
	The group $\Gamma$ contains a parabolic element 
	$\left[
			  \begin{matrix}
				1 & \varepsilon (\tr(u) + \tr(v))/c \\
				0 & 1
				\end{matrix}
	 \right]$. 
	 	
	\begin{figure}
  \centering
  \includegraphics[width=1.0\textwidth]{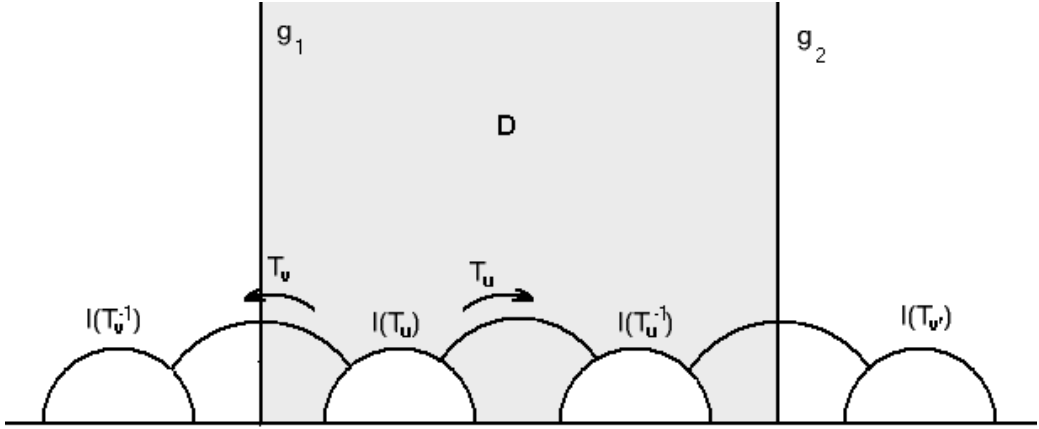}
   \caption{The case $\frac{a_1+d}{c}>0$.}
   \label{F:generationParab}
  \end{figure}
	
	Indeed,
	\[
	T:=T_{u} T_{v}^{-1} = \left[
			   \begin{matrix}
			    a_1 & b_1 \\
			    c & d
			   \end{matrix}
			   \right] 
			   \left[
			   \begin{matrix}
			    d & - b_2 \\
			    - c & a_2
			   \end{matrix}
			   \right]
			 = \left[
			   \begin{matrix}
			    1 &  b_1 a_2 - a_1 b_2 \\
			    0 & 1
			   \end{matrix}
			   \right]
	\]
	and
	\[
	(b_1 a_2 - a_1 b_2) c = (a_1 d -1) a_2 - a_1 (a_2 d - 1) = 
				  (a_1 + d) - (a_2 + d) = 
				 \varepsilon ( \tr(u) + \tr(v)).
	\]
	From $c \neq 0$ it follows that $I(T_u)$ and $I(T_v)$ exist and coincide.
	
	We notice that since $u$ and $v$ are nonnegative numbers then $\tr(u) \geq 2$ and $\tr(v) \geq 2$ and hence $\tr(u)$ and $\tr(v)$ cannot be elliptic transformations. Now we consider a region $D$ like the one indicated in Figure \ref{F:generationParab}. It is determined by the isometric circles of $T_{u}$ and $T_{u}^{-1}$ and by the geodesics $g_1$ and $g_2:=T(g_1)$, where
	\begin{itemize}
	 \item $g_1$ is the geodesic through $\infty$ orthogonal to $C(T_v)$, if $T_v$ is a hyperbolic transformation,
	 \item $g_1$ is the geodesic through $\infty$ and the fixed point of $T_v$, if $T_v$ is a parabolic transformation.
	\end{itemize}
	 From Poincare's theorem for fundamental polygons (see e.g. \cite{bM71}) it follows that $D$ is a fundamental domain for the Fuchsian group $\left\langle T,T_u\right\rangle=\left\langle T_u,T_v\right\rangle$.
	 We are going to show that $\Gamma \backslash D$ contains $Y(u,v,0)$.
	 \begin{itemize}
	 \item If $T_u$ is a hyperbolic transformation, then $C(T_u)$ is orthogonal to $I(T_u)$ and $I(T_u^{-1})$ (Lemma \ref{L:CIorthogonal}). Moreover $I(T_u^{-1})=T_u(I(T_u))$ and thus $\langle T_u \rangle  \backslash C(T_u)$ is a simple closed geodesic on $\Gamma\backslash \bbH$ with length $u$. 
	 \item If $T_u$ is a parabolic transformation, then the fixed point of $T_u$ corresponds to a cusp of $\Gamma \backslash \bbH$, i.e to a closed geodesic of length 0.
	 \item If $T_v$ is a hyperbolic transformation, then for $T_{v'}:=T_uT^{-1}$ we have $T_{v'} = TT_vT^{-1}$ and therefore $C(T_{v'}) = T(C(T_v))$, $I(T_{v'}) = T(I(T_v))$ and $I(T_{v'}^{-1}) = T(I(T_v^{-1}))$. Since $g_1$ is orthogonal to $C(T_v)$, $g_2= T(g_1)$ is orthogonal to $C(T_{v'})$ . Thus the part of $C(T_{v'})$ between $g_2$ and $I(T_{v'})$ (outside of $D$) is equivalent under $\Gamma$ to the part of $C(T_v)$ in $D$. Similarly the part of $C(T_v)$ between $g_1$ and $C(T_v^{-1})$ (outside of $D$) is equivalent to the part of $C(T_{v'})$ in $D$. This together with Lemma \ref{L:CIorthogonal} shows that $\langle T_u \rangle  \backslash C(T_v)$ is a simple closed geodesic on $\Gamma \backslash \bbH$ with length $v$. 
	 \item If $T_v$ is a parabolic transformation, then for $T_{v'}:=T_uT^{-1}$ we have $T_{v'} = TT_vT^{-1}$ and therefore the fixed point of $T_v$ is mapped to the fixed point of $T_{v'}$ by $T$ and thus both fixed points correspond to a cusp on $\Gamma \backslash \bbH$.
	 \end{itemize}
	 	\begin{figure}
  \centering
  \includegraphics[width=1.0\textwidth]{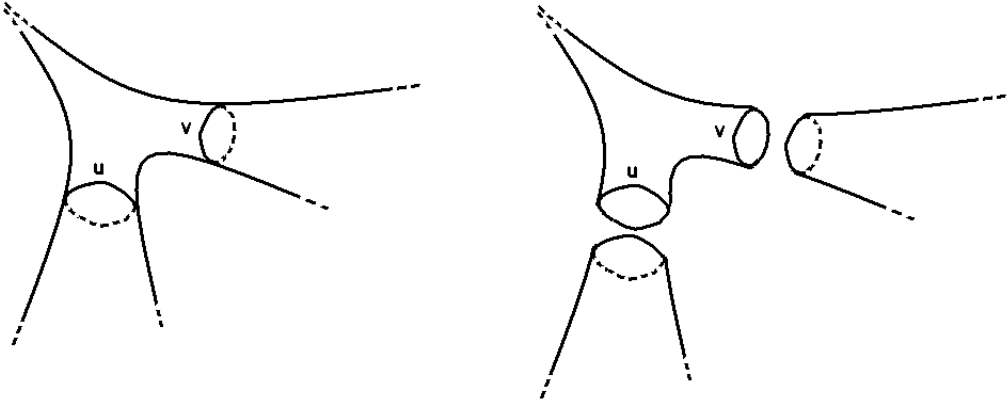}
   \caption{The surface $\Gamma \backslash \bbH$ and the Y-piece $Y(u,v,0)$.}
   \label{F:generationParab2}
  \end{figure}
	 Since $\left\langle T_u,T_v\right\rangle$ is a non-elementary Fuchsian group, the fixed points of $T_u$ and $T_v$ do not coincide. If both $T_u$ and $T_v$ are hyperbolic isometries, it follows from Lemma \ref{L:Notintersect} with $\lambda = 1$ that $C(T_u)$ and $C(T_v)$ do not intersect. Cutting $\Gamma \backslash \bbH$ along the closed geodesics described above of length $u$ and $v$, respectively, produces the required Y-piece $Y(u,v,0)$ (see Fig. \ref{F:generationParab2}).
	\end{proof}
  \begin{rem} The proposition remains true if $T_u$ is an elliptic transformation of finite order, i.e. $a_1+d \in (-2,2)$. Then the Fuchsian group $\left\langle T_u, T_v\right\rangle$ contains a degenerated Y-piece $Y(u,v,0)$, where $u$ is an elliptic fixed point.
  \end{rem}
	In the next corollary we use the notation of Proposition \ref{P:generationParab}.
	\begin{corollary}
	\label{C:1101}
	Let $\Gamma$ be a Fuchsian group containing the parabolic element $T = \begin{bmatrix} 1 & 1 \\ 0 & 1 \end{bmatrix}$. Then for every element $T_u = \begin{bmatrix} a_1 & b_1 \\ c & d \end{bmatrix}$, $c \neq 0$, there exists $T_v \in \Gamma$, $v\geq 0$, such that $\left\langle T_u, T_v\right\rangle \backslash \bbH$ contains $Y(u,v,0)$ with $\tr(u)= \left| a_1+d \right|$.
	\end{corollary}
	\begin{proof} For any $k \in \Z$ we consider
	$$ T^kT_u =  \begin{bmatrix} 1 & k \\ 0 & 1 \end{bmatrix}
	             \begin{bmatrix} a_1 & b_1 \\ c & d \end{bmatrix} =
	   \begin{bmatrix} a_1+kc & b_1+kd \\ c & d \end{bmatrix}.$$
	Pick $k' \in \Z$ such that $(a_1+d)(a_1+k'c +d) \leq 0$ and $|a_1+k'c +d| \geq 2$. Then set $T_v:=T^{k'}T_u$ and the claim follows from Proposition \ref{P:generationParab} and the previous remark.
	\end{proof}
	
	\begin{corollary}
	\label{C:AllInYpiece}
	Let $\Gamma$ be a Fuchsian group containing at least one parabolic element. Then for every non-parabolic element $T_u$ in $\Gamma$ there exists $T_v \in \Gamma$ such that $\left\langle T_u, T_v\right\rangle \backslash \bbH$ contains an Y-piece $Y(u,v,0)$ with $\tr(u) = \tr(T_u)$.
	\end{corollary}
  \begin{proof} If $\Gamma$ contains a parabolic element $T_1$ then, for some $R \in PSL(2, \R)$, $RT_1R^{-1} = T$ or $RT_1^{-1}R^{-1} = T$ where $T = \begin{bmatrix} 1 & 1 \\ 0 & 1 \end{bmatrix}$. If $R\Gamma R^{-1}$ contains also an element $A = \begin{bmatrix} a & b \\ 0 & d \end{bmatrix}$ then $A$ is a parabolic element because otherwise the group $\left\langle T,A\right\rangle$ would not be discrete. From Corollary \ref{C:1101} it follows that for every non-parabolic element $R T_u R^{-1}$ in $R\Gamma R^{-1}$ there exists $R T_v R^{-1} \in R \Gamma R^{-1}$ such that $\left\langle R T_u R^{-1}, R T_v R^{-1}\right\rangle \backslash \bbH$ contains an Y-piece $Y(u,v,0)$ with $\tr(u) = \tr(T_u)$. And hence for every non-parabolic element $T_u$ in $\Gamma$ there exists $T_v \in \Gamma$ such that $\left\langle T_u, T_v\right\rangle \backslash \bbH$ contains an Y-piece $Y(u,v,0)$ with $\tr(u) = \tr(T_u)$.
  \end{proof}
   
\subsection{Some geodesics on Y-pieces}
In this section we discuss several technical lemmas due to Schmutz  which we need in  the proof of Sarnak's conjecture.

	\begin{lemma}[ \cite{pS96}] \label{L:1}
	For all positive integers n, $Y(x, y, 0)$ contains $Y(\nu_{n}, y, 0)$, where
	\[
		\tr(\nu_{n}) = n ( \tr(x) + \tr(y) ) - \tr(y).
	\]
	In particular,  $\Tr(Y(x,y,0))$ contains the set $ \lbrace \tr(\nu_{n}):~ n = 1$, $2$, $3 \ldots \rbrace $.
	\end{lemma}
		
	For completeness we reproduce the
  proof given in \cite{pS96}.

	\begin{proof} 
\begin{figure}
\centering
\includegraphics[width=1.0\textwidth]{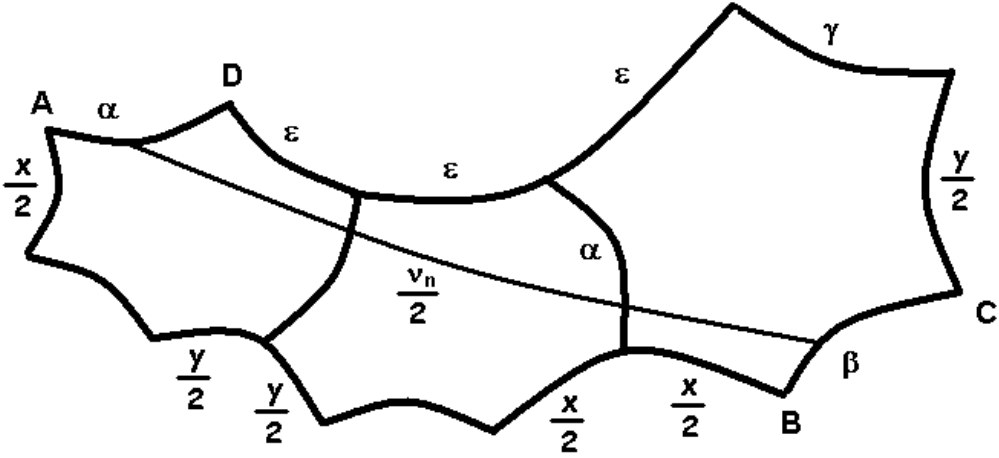}
\caption{The case $n=3$.}
\label{F:L1}
\end{figure}

	We replace $Y(x, y, 0)$ by $Y(x, y, 2\varepsilon)$, and we work on a covering surface of $Y(x, y, 2\varepsilon)$. If $x = 0$ or $y = 0$ we use the same trick.
	We consider  half of an $n$-fold covering $P$ of $Y(x,y,2\varepsilon)$ which is a convex geodesic $4 + 2 n$-gon (see Fig. \ref{F:L1}).  This implies that every quadrilateral ABCD (with sides AD and BC being also sides of $P$) lies in $P$ and from that we conclude that all angles of ABCD are smaller than or equal to $\frac{\pi}{2}$. Hence AD and BC have a common orthogonal $\frac{\nu_{n}}{2}$ which lies within ABCD. We will use the following formula (see for example \cite{aB95}): For any convex right-angled geodesic hexagon with consecutive sides $a_1$, $a_2$, $a_3$, $a_4$, $a_5$ and  $a_6$ holds:
	\[
		\cosh a_5 = \cosh a_2 \sinh a_1 \sinh a_3 - \cosh a_1 \cosh a_3
	\]
	
	In our case, we obtain on the one hand from the right-angled hexagon with consecutive sides $\frac{y}{2}$, $\gamma$, $\varepsilon$, $\alpha$, $\frac{x}{2}$, $\beta$ the equality
	\[
		\cosh \frac{x}{2} =
		  \cosh \gamma \sinh \frac{y}{2} \sinh \varepsilon -
		  \cosh \frac{y}{2} \cosh \varepsilon.
	\]
	On the other hand  the right-angled hexagon with consecutive sides $\frac{y}{2}$, $\gamma$, $ n \varepsilon$, a part of $\alpha$, $\frac{\nu_{n}}{2}$, a part of $\beta$, yields
	\[
		\cosh \frac{\nu_n}{2} =
		  \cosh \gamma \sinh \frac{y}{2} \sinh n \varepsilon - 
		  \cosh \frac{y}{2} \cosh n \varepsilon.
	\]
	Thus
	\[
		\cosh \frac{\nu_n}{2} =  
		  \frac {\cosh \frac{x}{2} + \cosh \frac{y}{2} \cosh \varepsilon}
		        {\sinh \frac{y}{2} \sinh \varepsilon}
		  \sinh \frac{y}{2} \sinh n \varepsilon - 
		  \cosh \frac{y}{2} \cosh n \varepsilon.
	\]
	Since $\lim_{\varepsilon \to 0} \frac{\sinh n \varepsilon}{\sinh \varepsilon} = n$, 
	\[
		\lim_{\varepsilon \to 0}\cosh \frac{\nu_n}{2} =  
		   n (\cosh \frac{x}{2} + \cosh \frac{y}{2}) - 
		  \cosh \frac{y}{2}.
	\]
	Thus the second right-angled hexagon determines in the limit case an Y-piece $Y(\nu_{n}, b, 0)$ with
	\[
		\tr(\nu_{n}) = n ( \tr(x) + \tr(y) ) - \tr(y).
	\]
	
	\end{proof}
	
	\begin{rem}
	Lemma \ref{L:1} is true even if $Y(x,y,0)$ is a degenerated Y-piece where $x$ corresponds to an elliptic fixed point and $y$ is a closed geodesic. Then $\tr(x)$ is equal to the trace of the generating elliptic ellement. For the proof we use again the half of an $n$-fold covering of $Y(x,y,2\varepsilon)$ which in this case  is a convex geodesic $4 + n$-gon and instead of the formula for a right-angled geodesic hexagon we use a similar formula for a geodesic pentagon with four right angles (\cite{aB95}, Theorem 7.18.1).
	\end{rem}
	
	The next three lemmas can be proved using  ideas simmilar to those in the above proof.
	
	\begin{lemma}[\cite{pS96}]
	\label{L:Y(yk,ym,0)}
	$\Tr(Y(x,0,0))$ contains $\Tr(Y(\lambda_k,\mu_m,0))$ with $\tr(\lambda_k)=k(\tr(x)+2)+2$ and $\tr(\mu_m)=m(\tr(x)+2)-2$ for all pairs $(k,m)$, $k,m \in\Z^+$.
	\end{lemma}
		
	\begin{lemma}[\cite{pS96}]
	\label{L:Y(a,2b,z)}
	$\Tr(Y(x,y,0))$ contains $\Tr(Y(\nu,2y,0))$ with $\tr(\nu)=2 + \tr(x)\tr(y).$
	\end{lemma}
			
	\begin{lemma}[\cite{pS96}]
	\label{L:Y(a,b,c)containsY(x,c,c)}
	$\Tr(Y(x,y,0))$ contains $\Tr(Y(\nu,0,0))$ where $$\tr(\nu) = (\tr(x)+\tr(y))^2 - 2.$$
	\end{lemma}
	

\section{The growth of the length spectrum}
\label{growth}
\setcounter{theorem}{0}

In  \cite{pS96} Schmutz proposes a proof of Conjecture \ref{Con:Schmutz} under the assumption that  the group $\Gamma$ contains parabolic elements. Unfortunately the proof contains 
a gap as we will explain in this section. However, using ideas and methods similar to those  in \cite{pS96} we are able to 
prove (part of) Sarnak's conjecture: Let $\Gamma$ be a cofinite Fuchsian group, which containins  parabolic elements.
If $\Gamma$ satisfies the B-C property, then $\Gamma$ is
 arithmetic. 

\subsection{An attempt to prove Conjecture \ref{Con:Schmutz}}
\label{SS:ProofAttempt}

We will need the following three results. 
Theorem 4.3 below is a direct corollary of Theorem \ref{T:TakeuchiArithmQuat}. Proofs of Theorem \ref{T:DivisionIsCompact} and Theorem \ref{T:DivisionAlgebra} can be found in \cite{sK92}, Chapter 5.

\begin{theorem}
\label{T:DivisionIsCompact}
Let $\Gamma$ be a Fuchsian group derived from a division quaternion algebra. Then the quotient space $\Gamma \backslash \bbH$ is compact, i.e. $\Gamma$ contains no parabolic elements.
\end{theorem}

\begin{theorem}
\label{T:DivisionAlgebra}
Let $A$ be a quaternion algebra over a totally real algebraic number field $F$ like in Theorem \ref{T:G(A,O)FuchsianGroup}. If $F \neq \Q$ then $A$ is a division quaternion algebra.
\end{theorem}

\begin{theorem}
\label{T:ArithmQuatOverQ}
Let $\Gamma$ be a cofinite Fuchsian group. Then $\Gamma$ is derived from a quaternion algebra over $\Q$ if and only if for every $T \in \Gamma$, $\tr(T) \in \Z$, i.e. $\Tr(\Gamma) \subseteq \Z$.
\end{theorem}

Now let $\Gamma$ be a cofinite Fuchsian group with at least one parabolic element.
If $\Gamma$ is derived from a quaternion algebra $A$, then $A$ is not a division quaternion algebra (Theorem \ref{T:DivisionIsCompact}) and consequently $A$ is a quaternion algebra over $\Q$ (Theorem \ref{T:DivisionAlgebra}). Hence if we prove that $\Gamma$ is derived from a quaternion algebra, then $\Gamma$ will be derived from a quaternion algebra over $\Q$. Hence, by Theorem \ref{T:ArithmQuatOverQ}, in order to prove the second part of Conjecture \ref{Con:Sarnak} in the case when $\Gamma$ contains at least one parabolic element it is enough to show that $Gap(\Gamma) > 0$ implies $\Tr(\Gamma) \subseteq \Z$.

If one wishes to show that $\Gamma$ is an arithmetic Fuchsian group it is enough to show that $\Gamma^{(2)}$ is derived from a quaternion algebra (Theorem \ref{T:TakeuchiArithm}). And since $\Gamma^{(2)}$ also contains at least one parabolic element it is sufficient to show that $\Tr(\Gamma^{(2)}) \subseteq \Z$ which is the same as to show that $\{\tr(a)^2 \mid a \in \Gamma \} \subseteq \Z$ because $\tr(a^2) = \tr(a)^2 - 2$.

The idea of a possible proof of Conjecture \ref{Con:Schmutz} is now the following:

\paragraph{ }
For an Y-piece $Y(a,b,c)$ we set $Gap(Y(a,b,c)):=Gap(\left\langle T_a,T_b\right\rangle)$, where $T_a$ and $T_b$ are isometries generating $Y(a,b,c)$ like in Proposition \ref{P:generationParab}.

\paragraph{ }
From Corollary \ref{C:AllInYpiece} we know that for every non-parabolic element $T_x$ in $\Gamma$ there exists $T_y \in \Gamma$ such that $\left\langle T_x, T_y\right\rangle \backslash \bbH$ contains an Y-piece $Y(x,y,0)$ with $\tr(x) = \tr(T_x)$. Since the trace of every parabolic transformation is equal to 2, it is enough to show that if $\Tr(\Gamma)$ has linear growth then, for every $Y(x,y,0)$, $\tr(x)^2$ and $\tr(y)^2$ are integers. 

In \cite{pS96} Schmutz proves the following two propositions:

\begin{proposition}[\cite{pS96}]
\label{P:GapTr(x)Integer}
$Gap(Y(x,0,0))>0$ if and only if $\tr(x)$ is an integer.
\end{proposition}

\begin{proposition}[\cite{pS96}]
\label{P:GapSquaresIntegers}
If $Gap(Y(x,y,0))>0$ then the numbers $\tr(x)^2$, $\tr(y)^2$ and $\tr(x)\tr(y)$ are integers.
\end{proposition}

\begin{proof}[Idea of the proof:] From Lemma \ref{L:Y(a,b,c)containsY(x,c,c)} it follows that, for an Y-piece $Y(a,b,0)$, $\Tr(Y(a,b,0))$ contains $\Tr(Y(z,0,0))$ with $\tr(z) = (\tr(a) + \tr(b))^2 -2$. If $Y(x,y,0)$ contains $Y(a,b,0)$ then $Y(x,y,0)$ contains also $Y(z,0,0)$. Hence $Gap(Y(z,0,0))>0$ and by Proposition \ref{P:GapTr(x)Integer} $\tr(z)$ is an integer. 
The proposition is then proved by applying the above considerations to different Y-pieces $Y(a,b,0)$ contained in $Y(x,y,0)$.
\end{proof}

Observe that the condition $Gap(Y(x,y,0))>0$ is used only in case we need $Gap(Y(z,0,0))>0$ in order to apply Proposition \ref{P:GapTr(x)Integer} for an Y-piece $Y(z,0,0)$ contained in $Y(x,y,0)$.

If $\Tr(\Gamma)$ has linear growth then $\Tr(Y(x,y,0))$ has linear growth for every Y-piece $Y(x,y,0)$ contained in $\Gamma \backslash \bbH$. Our aim is to prove that if $\Tr(Y(x,y,0))$ has linear growth then $\tr(x)^2$ and $\tr(y)^2$ are integers.

The idea of Schmutz is to proceed as in  the proof of Proposition \ref{P:GapSquaresIntegers}, but instead of
Proposition \ref{P:GapTr(x)Integer} to use the following

\begin{claim}
\label{Cl:LinearGrowthTr(x)Integer}
$\Tr(Y(x,0,0))$ has linear growth if and only if $\tr(x)$ is an integer.
\end{claim}

 Proposition \ref{P:GapTr(x)Integer} shows that  $Gap(Y(x, 0, 0)) >0$  if $\tr(x)$ is an integer and hence $\Tr(Y(x,0,0))$ has linear growth. So in order to prove Claim \ref{Cl:LinearGrowthTr(x)Integer} it remains to show that if $\Tr(Y(x,0,0))$ has linear growth then $\tr(x) \in \N$, which is the same as to show that if $\tr(x)$ is not an integer then $\Tr(Y(x,0,0))$ has not linear growth. 

If the real number $\tr(x)$ is not an integer it can be either rational or irrational. In the next two subsections we are going to present the proof of Claim \ref{Cl:LinearGrowthTr(x)Integer}  in \cite{pS96} in the case when $\tr(x)$ is not rational and to show that there is a gap in the proof of Claim \ref{Cl:LinearGrowthTr(x)Integer}  in \cite{pS96} in the case when $\tr(x)$ is rational.

\subsection{The proof of Claim \ref{Cl:LinearGrowthTr(x)Integer}  in \cite{pS96} in the case when $\tr(x)$ is not rational}
\label{S:NotRational}
\begin{proof}[We give the details of the proof in \cite{pS96}:]
We assume that $z:=\tr(x) + 2$ is not rational. By Lemma \ref{L:Y(yk,ym,0)}, $\Tr(Y(x,0,0))$ contains $\Tr(Y(\lambda_k,\mu_m,0))$ with $\tr(\lambda_k)=k(\tr(x)+2)+2$ and $\tr(\mu_m)=m(\tr(x)+2)-2$ for all pairs $(k,m)$, $k,m \in\Z^+$. Hence  it follows from Lemma \ref{L:Y(a,2b,z)} that $\Tr(Y(x,0,0))$ contains $\tr(\mu_m)\tr(\lambda_k)+2$ and thus the set
\[ \{mkz^2 - 2(k-m)z -2 \mid m,k \in \Z^+ \}. \]
We claim that for all different pairs of positive integers $(m_1,k_1)$ and $(m_2,k_2)$
\[ m_1k_1z^2 - 2(k_1 - m_1)z - 2 \neq m_2k_2z^2 - 2(k_2 - m_2)z - 2.\]
To see this we assume that
\[ m_1k_1z^2 - 2(k_1 - m_1)z - 2 = m_2k_2z^2 - 2(k_2 - m_2)z - 2.\]
Since $z \notin \Q$, we have $z \neq 0$ and the above equality is equivalent to
\[ (m_1k_1 - m_2k_2)z - 2(k_1 - m_1 - (k_2 - m_2)) = 0.\]
Now, if $Az + B = 0$ for some integers $A$ and $B$ and $z \notin \Q$, then $A=0$ and thus $B=0$. In our case this means $m_1k_1 - m_2k_2 = 0$ and $k_1 -m_1 = k_2 - m_2$. Consequently $k_1$, $-m_1$ and $k_2$, $-m_2$ are solutions of the  quadratic equation
\[ \alpha^2 - (k_1 - m_1)\alpha + m_1k_1 = 0.\]
Since $m_1$, $k_1$, $m_2$ and $k_2$ are positive, it follows that $k_1 = k_2=:k$ and $m_1 = m_2=:m$, a contradiction.

We assume that $k \geq m$. Then $mkz^2 \geq mkz^2 - 2(k-m)z -2$. Every $i \in \N$ can be written as a product of two positive integers in $\left[ \frac{1 + \sigma_0(i)}{2} \right]$ different ways, where $\sigma_0(i)$ is the number of different positive divisors of $i$. This implies that for each $N \in \N$
\begin{eqnarray*}
& & \#\{ a \in \Tr(Y(x,0,0)) \mid a \leq Nz^2 \} \\
& \geq & \#\{a := mkz^2 - 2(k-m)z -2 \mid a \leq Nz^2~\mathrm{and}~ k \geq m, ~k,m \in \Z^+ \} \\
& \geq & \#\{ mkz^2 \mid mkz^2 \leq Nz^2~\mathrm{and}~ k \geq m, ~k,m \in \Z^+\}\\
& \geq & \frac{1}{2} \sum^{N}_{i = 1} \sigma_0(i).
\end{eqnarray*}
If we can show that $\sum^{N}_{i = 1} \sigma_0(i)$ grows like $N\log N$, then $\Tr(Y(x,0,0))$ does not have linear growth.

In the sum $\sum^{N}_{i = 1} \sigma_0(i)$, $1$ is counted as a divisor $N$ times, 2 is counted $\left[ \frac{N}{2} \right]$ times, every integer $j \leq N$ is counted $\left[ \frac{N}{j} \right]$ times and therefore
\[ \sum^{N}_{i = 1} \sigma_0(i) = \sum^{N}_{j = 1} \left[ \frac{N}{j} \right]. \]
As $\frac{N}{j}-1 \leq \left[ \frac{N}{j} \right] \leq \frac{N}{j}$ and hence we have
\[ \sum^{N}_{j = 1} \frac{N}{j} - N \leq \sum^{N}_{j = 1} \left[ \frac{N}{j} \right] \leq \sum^{N}_{j = 1} \frac{N}{j}.\]
With the lower and upper Darboux sums for the function $f(x) = \frac{1}{x}$ in the interval $[1, N]$ with the partition of the interval given by the integers between $1$ and $N$, we obtain the following inequalities:
\[ \sum^{N}_{j = 2} \frac{1}{j} \leq \int^{N}_{1} \frac{1}{x}dx \leq \sum^{N-1}_{j = 1}\frac{1}{j}.\]
Since $\int^{N}_{1} \frac{1}{x}dx = \log N - \log 1 = \log N$, we have
\[ \sum^{N}_{j = 1} \frac{N}{j} \geq N\sum^{N-1}_{j = 1}\frac{1}{j} \geq N \log N ~~~~~ \mathrm{and}\]
\[ N(\log N + 1) \geq N(\sum^{N}_{j = 2} \frac{1}{j} + 1) = \sum^{N}_{j = 1} \frac{N}{j}. \]
Hence
\[ N\log N + N \geq \sum^{N}_{i = 1} \sigma_0(i) = \sum^{N}_{j = 1} \left[ \frac{N}{j} \right] \geq N\log N - N \]
which means that $\sum^{N}_{i = 1} \sigma_0(i)$ grows like $N \log N$ and in particular not linear (and does not satisfy the B-C property). This proves 
Claim 4.6  in the case when $z$ is not a rational number.
\end{proof}

\subsection{The gap in the proof of Claim \ref{Cl:LinearGrowthTr(x)Integer} in \cite{pS96} in the case when $\tr(x)$ is rational}
\label{S:Rational}
Unfortunately we cannot use the above argument in the case when $z=tr(x)+2$ is a rational number $\frac{a}{b}$ with $b>1$ and $(a,b)=1$, because $v_1:=m_1k_1z^2 - 2(k_1 - m_1)z - 2$ and  $v_2 := m_2k_2z^2 - 2(k_2 - m_2)z - 2$ can be equal for different pairs $(k_1,m_1)$ and $(k_2,m_2)$. 

Indeed,  assume that $v_1 = v_2$ or equivalently, since $z >0$,
\[ (m_1k_1 - m_2k_2)z -2(k_1 - m_1 - (k_2 - m_2)) = 0. \]
Now, if $Az + B = 0$ for some integers $A$ and $B$ and $z = \frac{a}{b}$, then $A$ must not be $0$, it can also be divisible by $b$. But if $|A|<b$ then $A=0$ and thus $B=0$ and as in Section  \ref{S:NotRational} we have $k_1 = k_2$ and $m_1 = m_2$.

Therefore, since $k_1$, $m_1$, $k_2$ and $m_2$ are positive, we can guarantee that $v_1$ and $v_2$ are different for different pairs $(k_1,m_1)$ and $(k_2,m_2)$, if $m_1k_1 < b$ and $m_2k_2 <b$ and thus as in Section \ref{S:NotRational} we get
\[ \#\{ y \in \Tr(Y(x,0,0)) \mid y \leq bz^2 \} \geq \frac{1}{2}\sum^{b}_{i=1} \sigma_0(i) \geq \frac{1}{2}(b\log b - b). \]
From Lemma \ref{L:Y(a,b,c)containsY(x,c,c)} it follows that $\Tr(Y(x,0,0))$ contains $\Tr(Y(x_2,0,0))$ with $\tr(x_2) = (\tr(x)+2)^2 -2 = z^2 - 2$. By induction $\Tr(Y(x,0,0))$ contains $\Tr(Y(x_n,0,0))$ with $\tr(x_n) = z^{(2^n)}-2$.

In \cite{pS96} the author suggests to use the above estimates of the trace set for every $Y(x_n,0,0)$ (in this case $\tr(x_n)+2 = z^{2^n} = \frac{a^{2^n}}{b^{2^n}}$):
\[ \#\{ y \in \Tr(Y(x_n,0,0)) \mid y \leq b^{2^n}z^{2^{n+1}} \} \geq \frac{1}{2}\sum^{b^{2^n}}_{i=1} \sigma_0(i) \geq \frac{1}{2}b^{2^n}(\log b^{2^n} - b^{2^n}).\]
He claims that $\Tr(Y(x,0,0))$ has not linear growth because for every $n \in \N$
\begin{eqnarray*}\# \{y \in \Tr(Y(x,0,0)) \mid y \leq b^{2^n}z^{2^{n+1}}\} & \geq & \# \{y \in \Tr(Y(x_n,0,0)) \mid y \leq b^{2^n}z^{2^{n+1}}\}\\
& \geq & \frac{1}{2}\sum^{b^{2^n}}_{i=1} \sigma_0(i) \geq \frac{1}{2}b^{2^n}(\log b^{2^n} - b^{2^n}).
\end{eqnarray*} 
If $z^{2^{n+1}}$ were a constant then this argumentation would work. However $z^{2^{n+1}}$ also grows when $n$ grows.

An immediate counter-example are the Y-pieces $Y(x=z-2,0,0)$ with $z^2 = \frac{a^2}{b^2} > b$: If the estimate
\[ \# \{y \in \Tr(Y(x,0,0)) \mid y \leq b^{2^n}z^{2^{n+1}}\} \geq \frac{1}{2}\sum^{b^{2^n}}_{i=1} \sigma_0(i)
\]
implies non-linear growth then there exists $n_0 \in \N$ such that for infinitely many $n \geq n_0$ the inequality $b^{2^n}z^{2^{n+1}} \leq \frac{1}{2}\sum^{b^{2^n}}_{i=1} \sigma_0(i)$ holds. But this is not the case when $z^2 > b$. In fact, for all positive integers $n$, one has in that case
\[ b^{2^n}z^{2^{n+1}} = b^{2^n} (z^2)^{2^n} > b^{2^n} b^{2^n} > \frac{1}{2}b^{2^n} (\log b^{2^n}+1) \geq \frac{1}{2}\sum^{b^{2^n}}_{i=1} \sigma_0(i).\]

At first view a possible reason why the above considerations did not suffice to prove the non-linear growth of $\Tr(Y(x,0,0))$ might be that 
 not enough  elements of the set 
\[ S_n = \{ mkz^{2^{n+1}} - 2(k-m)z^{2^n} -2 \mid m,k \in \Z^+ \} \]
have been taken into accout.
But it turns out that even in the union $\bigcup_{n=0}^{\infty}S_n$
there are not enough different numbers to guarantee  non-linear growth of  $\Tr(Y(x,0,0))$. Indeed, every $y \in S_0$ has the form 
\[ mk\frac{a^2}{b^2} - 2(k-m)\frac{a}{b} -2 = \frac{a}{b^2}\left(mka - 2(k-m)b\right)-2. \]
Hence
\[ S_0 \subseteq B_0 := \{v:= \frac{a}{b^2}j -2 \mid j \in \N,~ v > 0 \}.\]
The number of the elements in $B_0$ which are smaller than $N \in \N$ is bounded by $\frac{N+2}{\frac{a}{b^2}} = (N+2) \frac{b^2}{a}$.

Analogously we get for every $n \in \N$ and $N \in \N$ and $B_n = \{v:=\left(\frac{a}{b^2}\right)^{2^n}j -2 \mid j \in \N,~ v > 0 \}$
\[ \#\{w \in S_n \mid w \leq N\} \leq \#\{v \in B_n \mid v \leq N\} \leq (N+2) \left(\frac{b^2}{a}\right)^{2^n}. \]
Hence
\[ \#\{w \in \bigcup^{\infty}_{n=0}S_n \mid w \leq N\} \leq \#\{v \in \bigcup^{\infty}_{n=0}B_n \mid v \leq N\} \leq (N+2) \sum^{\infty}_{n=0}\left(\frac{b^2}{a}\right)^{2^n}.\]
If $a > b^2$ the last sum is convergent and independent of $N$, i.e.
\[ \#\{w \in \bigcup^{\infty}_{n=0}S_n \mid w \leq N\} \leq const (N+2) \]
which means that $\bigcup^{\infty}_{n=0}S_n$ has only linear growth! 
Thus if $\tr(x)$ is rational the previous argument due to Schmutz is not conclusive: $\tr(x)\in \Q \backslash \Z$ does \emph{not} necessarily imply that $\Tr(Y(x,0,0))$ does not grow linearly! However, we will see in the next section that $\tr(x)\in \Q \backslash \Z$ implies that $\Tr(Y(x,0,0))$ does not satisfy the B-C property.

\subsection{$\Tr(Y(x,0,0))$ satisfies the B-C property if and only if $\tr(x)$ is an integer}

The aim of this section is the proof of the following Theorem \ref{T:Sarnak1}, which in turn
 proves the first part of Conjecture \ref{Con:Sarnak}.
\begin{theorem}
\label{T:Sarnak1}
Let $\Gamma$ be a cofinite Fuchsian group with at least one parabolic element. Then $\Tr(\Gamma)$ satisfies the B-C property if and only if $\Gamma$ is an arithmetic group.
\end{theorem}
\begin{proof}
By Theorem \ref{T:LuoSarnak} if $\Gamma$ is an arithmetic group then $\Tr(\Gamma)$ satisfies the B-C property. So it remains to show that if $\Tr(\Gamma)$ satisfies the B-C property then $\Gamma$ is an arithmetic Fuchsian group. The proof below follows the ideas of Section \ref{SS:ProofAttempt} but instead of Claim \ref{Cl:LinearGrowthTr(x)Integer} we use Proposition \ref{P:B-CPropertyTr(x)Integer}.

By Corollary \ref{C:AllInYpiece}, for every non-parabolic element $T_x$ in $\Gamma$ there exists $T_y \in \Gamma$ such that $\left\langle T_x, T_y\right\rangle \backslash \bbH$ contains an Y-piece $Y(x,y,0)$ with $\tr(x) = \tr(T_x)$. If $T_x$ is an elliptic element, then $Y(x,y,0)$ is a degenerated Y-piece. By \S \ref{SS:ProofAttempt} it is enough to show that if $\Tr(\Gamma)$ satisfies the B-C property then, for every $Y(x,y,0)$, $\tr(x)^2$ and $\tr(y)^2$ are integers.

If $\Gamma$ satisfies the B-C property then, for every Y-piece $Y(x,y,0)$ contained in $\Gamma \backslash \bbH$, the trace set $\Tr(Y(x,y,0))$ also satisfies the B-C property. Hence it is enough to show that if $\Tr(Y(x,y,0))$ satisfies the B-C property then $\tr(x)^2$ and $\tr(y)^2$ are integers. 

If $Y(x,y,0)$ is non-degenerated then the claim follows from the next Proposition \ref{P:B-CPropertySquaresIntegers}.

If $Y(x,y,0)$ is degenerated, i.e. $x$ corresponds to an elliptic fixed point, then by the remark after Lemma \ref{L:1} the Y-piece $Y(x,y,0)$ contains $Y(\nu_2,y,0)$ and $Y(\nu_3,y,0)$ with $\tr(\nu_2)= 2 \tr(x) + \tr(y)$ and $\tr(\nu_3) = 3 \tr(x) + 2 \tr(y)$. Since $\tr(y) \geq 2$ then $\tr(\nu_2)$ and $\tr(\nu_3)$ are also greater or equal $2$. Hence $Y(\nu_2,y,0)$ and $Y(\nu_3,y,0)$ are non-degenerated and by the next Proposition \ref{P:B-CPropertySquaresIntegers} it follows that $\tr(\nu_2)^2$, $\tr(\nu_3)^2$ and $\tr(y)^2$ are integers.

So $4\tr(x)^2 + 4\tr(x)\tr(y) = \tr(\nu_2)^2 - \tr(y)^2$ and $3\tr(x)^2 + 4\tr(x)\tr(y) = \tr(\nu_3)^2 - \tr(y)^2$ are integers and hence $\tr(x)^2$ is an integer.
\end{proof}

\begin{proposition}
\label{P:B-CPropertySquaresIntegers}
If $\Tr(Y(x,y,0))$ satisfies the B-C property then $\tr(x)^2$, $\tr(y)^2$ and $\tr(x)\tr(y)$ are integers.
\end{proposition}

\begin{proof}
The proof is the same as that  of Proposition \ref{P:GapSquaresIntegers} but instead of Proposition \ref{P:GapTr(x)Integer} we use Proposition \ref{P:B-CPropertyTr(x)Integer} below.
\end{proof}

\begin{proposition}
\label{P:B-CPropertyTr(x)Integer}
$\Tr(Y(x,0,0))$ satisfies the B-C property if and only if $\tr(x)$ is an integer.
\end{proposition}

In the rest of this Section we are going to 
 prove Proposition \ref{P:B-CPropertyTr(x)Integer}. We will need the following Lemma:
\begin{lemma}
\label{L:Bezout}
Let $a$ and $b$ be coprime natural numbers, which are greater than $1$. Further let $b = p b_1$, where $p$ is a prime number and $b_1 \in \N$. Then there exist $u, v \in \N \backslash \{0\}$ such that $|ua - vb| = 1$, $v < a$ and $(v,p)=1$ (and thus also $(v,a)=1$ and $(u,b)=1$).
\end{lemma}
\begin{proof}
Bezout's identity yields $u', v' \in \Z \backslash \{0\}$ such that $u'a+v'b = 1$. We can also write this equivalently as $|\tilde{u}a - \tilde{v}b| = 1$, where  $\tilde{u}$ and $\tilde{v}$ are positive natural numbers. Furthermore, we have that $\tilde{v} = qa + r$, where $q,r \in \N$, $r < a$ and $r > 0$, because $(\tilde{v},a) = 1$ and $a > 1$. Thus after subtracting  $0 = q(ba - ab)$ from $|\tilde{u}a - \tilde{v}b|$ we get:
\[ |(\tilde{u} - qb)a - rb | = 1.\]

If $(r, p ) = 1$ we set $u:=\tilde{u} - qb$ and $v:=r$. Note that $u$ is positive because $a$ is positive and $rb > 1$.

If $(r, p ) = p$ then we subtract $0=ba - ab$ from $(\tilde{u} - qb)a - rb$. We obtain $|(\tilde{u} - (q+1)b)a + (a - r)b| =1$, where $0 < a - r < a$ and $(a - r, p) = 1$, because $(a,p) = 1$ (since $p$ is a divisor of $b$). From $(a - r)b > 1$ and $a>0$ it follows that  $\tilde{u} - (q+1)b < 0$. We set $ u = -(\tilde{u} - (q+1)b)$ and $v = a - r$.
\end{proof}

\begin{proof}[Proof of Proposition \ref{P:B-CPropertyTr(x)Integer}.]
If $\tr(x)$ is an integer then it follows from Proposition \ref{P:GapTr(x)Integer} that $Gap(Y(x,0,0))>0$. This means that in every interval $[n, n+1]$ there are at most $\left[\frac{1}{Gap(Y(x,0,0))} +1\right]$ elements from the set $\Tr(Y(x,0,0))$ and hence $\Tr(Y(x,0,0))$ satisfies the B-C property.

Now let $\Tr(Y(x,0,0))$ satisfy the B-C property. We assume that $\tr(x)$ is not an integer. There are two possibilities for $\tr(x)$:
\begin{itemize}
\item[\textbf{Case 1:}] $\tr(x)$ is not a rational number.
\end{itemize}
In section \ref{S:NotRational} we already showed that in this case $\Tr(Y(x,0,0))$ does not have linear growth and, in particular, does not satisfy the B-C property. Hence $\tr(x)$ cannot be irrational.
\begin{itemize}
\item[\textbf{Case 2:}] $\tr(x)$ is a rational number (but not an integer).
\end{itemize}
Then the number $z = \tr(x) + 2$ is equal to $\frac{a}{b}$ with $a$ and $b$ coprime natural numbers, $b > 1$ and $a > b$ because $z > 2$.

As in \S \ref{S:Rational} it follows from Lemma \ref{L:Y(a,b,c)containsY(x,c,c)} that $\Tr(Y(x,0,0))$ contains $\Tr(Y(x_k, 0, 0))$ with $\tr(x_k) = z^{2^k} -2$, $k \in \N$. By Lemma \ref{L:1} $\Tr(Y(x_k, 0, 0))$ contains the set
\[ \left\{m(z^{2^k} - 2 + 2) -2 \mid m \in \N \backslash \{0\} \right\} = \left\{m\left(\frac{a}{b}\right)^{2^k} -2 \mid m \in \N \backslash \{0\} \right\}.\]

We are going to show that for every $n \in \N$ there exist $n$ different numbers $z_{m_i,k_i}:= m_i\left(\frac{a}{b}\right)^{2^{k_i}}-2$, $i = 1,\ldots,n$, such that 
\[ \max \left\{z_{m_i,k_i} \mid i = 1, \ldots, n\right\}-\min \left\{z_{m_i,k_i} \mid i = 1, \ldots, n\right\} \leq 1.\]
And thus we show that $\Tr(Y(x,0,0))$ does not satisfy the B-C property.

\begin{itemize}
\item[\textit{Step 1.}] First we consider a function $f:\N \to \N$ with the following properties: $f(0) = 0$ and for $n>0$, $b^{2^{f(n)}} > 2 \prod^{n-1}_{i = 0}a^{2^{f(i)}}$.
\end{itemize}
Such function $f$ exists, because if we assume that we have defined $f$ for $0,\ldots,n-1$ then the right-hand side of the inequality is fixed and we can choose $f(n)$ big enough so that the inequality holds. We notice that  $f(n+1) > f(n)$ for every $n \in \N$ because
\[ b^{2^{f(n+1)}} > 2 \prod^{n}_{i = 0}a^{2^{f(i)}} \geq  a^{2^{f(n)}} > b^{2^{f(n)}}.\]

For convenience we set $g(n):= 2^{f(n)}$. Then we have $g(0)=1$ and for $n>0$, $b^{g(n)} > 2 \prod^{n-1}_{i = 0}a^{g(i)}$.
\begin{itemize}
\item[\textit{Step 2.}] We fix an arbitrary natural number $n$ greater than $1$. Let $b = pb_1$ where $p$ is a prime number and $b_1 \in \N$.
\item[\textit{Step 3.}] We can find positive integers $u_i,v_i$, $i = 1, \ldots, n$, such that
\[\left|u_i\left(\frac{a}{b}\right)^{g(i)} - v_i v_{i+1}\ldots v_n\frac{a}{b}\right| = \frac{a}{b^{g(i)}},\]
where $v_i < a^{g(i)-1}$, $(v_i,a)=1$ and $(v_i,p)=1$. In fact,
\end{itemize}
by Lemma \ref{L:Bezout} there exist $u_n, v_n \in \N \backslash \{0\}$ such that $\left|u_n a^{g(n)-1} - v_n b^{g(n)-1}\right|=1$, $v_n < a^{g(n)-1}$, $(v_n,a) = 1$ and $(v_n, p)=1$. Hence
\[ \left|u_n\left(\frac{a}{b}\right)^{g(n)} - v_n \frac{a}{b} \right| = \frac{a}{b^{g(n)}}\left|u_n a^{g(n)-1} - v_nb^{g(n)-1}\right| = \frac{a}{b^{g(n)}}.\]

Since $(a^{g(n-1)-1}, v_nb^{g(n-1)-1})=1$, then by Lemma \ref{L:Bezout} there exist $u_{n-1}, v_{n-1} \in \N \backslash \{ 0 \}$ such that $\left|u_{n-1} a^{g(n-1)-1} - v_{n-1}v_n b^{g(n-1)-1}\right|=1$, where $v_{n-1} < a^{g(n-1)-1}$, $(v_{n-1},a) = 1$ and $(v_{n-1}, p)=1$. Hence
\[ \left|u_{n-1}\left(\frac{a}{b}\right)^{g(n-1)} - v_{n-1}v_n \frac{a}{b} \right| = \frac{a}{b^{g(n-1)}}\left|u_{n-1} a^{g(n-1)-1} - v_{n-1}v_nb^{g(n-1)-1}\right| = \frac{a}{b^{g(n-1)}}. \]

For $1 \leq i \leq n-1$ we assume that $u_j, v_j$ are defined for all $j = i+1,\ldots,n$. We define $u_i$ and $v_i$:

Since $(a^{g(i)-1}, v_{i+1}\ldots v_{n-1}v_nb^{g(i)-1})=1$, then again by Lemma \ref{L:Bezout} there exist $u_i, v_i \in \N\backslash \{ 0 \}$ such that $\left|u_i a^{g(i)-1} - v_i v_{i+1}\ldots v_{n-1}v_nb^{g(i)-1}\right|$, where $v_i < a^{g(i)-1}$, $(v_i,a) = 1$ and $(v_i, p)=1$. Hence
\[ \left|u_i\left(\frac{a}{b}\right)^{g(i)} - v_i v_{i+1}\ldots v_{n-1}v_n \frac{a}{b} \right| = \frac{a}{b^{g(i)}}\left|u_i a^{g(i)-1} - v_i v_{i+1}\ldots v_{n-1}v_nb^{g(i)-1}\right| = \frac{a}{b^{g(i)}}. \]
\begin{itemize}
\item[\textit{Step 4.}] Set $m_0 := v_1\ldots v_{n-1}v_n$ and $m_i := v_1 \ldots v_{i-1}u_i$ for all $i = 1,\ldots,n$. We claim that the numbers $z_{m_i,f(i)} = m_i\left(\frac{a}{b}\right)^{2^{f(i)}}-2$, $i = 0, \ldots, n$, are all inside an interval of length $1$.
\end{itemize}
Indeed, for every $i = 1,\ldots,n$:
\begin{eqnarray*}
\left|z_{m_i,f(i)} - z_{m_0,f(0)}\right| & = & \left|m_i\left(\frac{a}{b}\right)^{g(i)}-2 - m_0\left(\frac{a}{b}\right)^{g(0)} +2\right|\\
& = & \left| v_1 \ldots v_{i-1}u_i \left(\frac{a}{b}\right)^{g(i)} - v_1 \ldots v_{i-1}v_i \ldots v_n \frac{a}{b} \right|\\
& = & v_1 \ldots v_{i-1} \left| u_i \left(\frac{a}{b}\right)^{g(i)} - v_i \ldots v_n \frac{a}{b} \right| = v_1 \ldots v_{i-1} \frac{a}{b^{g(i)}}\\
& < & \frac{a^{g(1)-1}\ldots a^{g(i-1)-1} a}{b^{g(i)}} \leq \frac{\prod^{i-1}_{j=0} a^{g(j)}}{b^{g(i)}} < \frac{1}{2},
\end{eqnarray*}
where the last inequality follows from our choice of the function $g$.
\begin{itemize}
\item[\textit{Step 5.}] We finally show that the numbers $z_{m_i,f(i)}$, $i = 0, \ldots, n$ are all different.
\end{itemize}
For every $i = 1, \ldots, n$, $u_i$ satisfies $(u_i,p)=1$, because $(u_i,b)=1$ (otherwise the difference $\left|u_ia^{g(i)-1}-v_i \ldots v_n b^{g(i)-1}\right|$ could not be equal to 1). We have chosen $v_i$, $i = 1, \ldots, n$, such that $(v_i,p)=1$. Hence $p$ does not divide $m_i = v_1 \ldots v_{i-1}u_i$ and $m_0 = v_1 \ldots v_{n-1}v_n$. Note also that $(p,a)=1$.

Let $d$ be the exponent of $p$ in the prime number decomposition of $b$. Write $z_{m_i,f(i)} = m_i\left(\frac{a}{b}\right)^{2^{f(i)}} - 2 = \frac{s}{t}$, with $s,t \in \N \backslash \{0\}$, $(s,t)=1$. Then $p^{d 2^{f(i)}}$ divides $t$ and $p^{d 2^{f(i)}+1}$ does not divide $t$.

Hence, since for $i \neq j$, $f(i) \neq f(j)$, the numbers $z_{m_i,f(i)}$ and $z_{m_j,f(j)}$ are different.
\end{proof}

\end{document}